\documentclass[11pt]{article}
\usepackage{times}
\usepackage{amsmath,amsfonts,amstext,latexsym,amssymb,amsbsy,amsopn,amsthm,eucal,slashed,amstext,enumerate,accents}
\usepackage{mathrsfs}
\usepackage[T1]{fontenc}

\usepackage{cases}
\usepackage{color}
\usepackage{txfonts}
\usepackage{dsfont}
\usepackage{graphicx}  
\usepackage{comment}

\allowdisplaybreaks[4]
\usepackage{scalerel,stackengine}
\stackMath
\newcommand\reallywidehat[1]{%
	\savestack{\tmpbox}{\stretchto{%
			\scaleto{%
				\scalerel*[\widthof{\ensuremath{#1}}]{\kern-.6pt\bigwedge\kern-.6pt}%
				{\rule[-\textheight/2]{1ex}{\textheight}}
			}{\textheight}%
		}{0.5ex}}%
	\stackon[1pt]{#1}{\tmpbox}%
}
\usepackage[colorlinks,
linkcolor=blue,      
anchorcolor=blue,  
citecolor=blue,        
]{hyperref}

\allowdisplaybreaks[4]

\begin{document}
	
	\setlength{\abovedisplayshortskip}{3mm}
	\setlength{\belowdisplayshortskip}{3mm}
	\setlength{\abovedisplayskip}{3mm}
	\setlength{\belowdisplayskip}{3mm}
	
	\newtheorem{theorem}{Theorem}[section]
	\newtheorem{corollary}[theorem]{Corollary}  
	\newtheorem{lemma}[theorem]{Lemma}
	\newtheorem{definition}[theorem]{Definition}
	\newtheorem{question}[theorem]{Question}
	\newtheorem{proposition}[theorem]{Proposition}
	\newtheorem{remark}[theorem]{Remark}
	\newtheorem{example}[theorem]{Example}
	\newtheorem{assumption}[theorem]{Assumption}
	\newtheorem{theoremalph}{Theorem}
	\renewcommand\thetheoremalph{\Alph{theoremalph}}

	\newenvironment{sequation}{\begin{equation}\small}{\end{equation}}
	\newenvironment{tequation}{\begin{equation}\tiny}{\end{equation}}
	\renewcommand{\thefootnote}{\fnsymbol{footnote}}


	\title
	{\bf\Large Fourth-Moment Theorems for Sums of Multiple Wigner Integrals}
	
	\author
	{\ \ Zizhan Chen,\quad 
		Xiang-Dong Li 
	}
	\date{}
	
	\maketitle
	\numberwithin{theorem}{section}
	\numberwithin{equation}{section}
	\setcounter{tocdepth}{2}
	\setcounter{secnumdepth}{2}
	
	\begin{abstract}
	We derive quantitative free fourth-moment bounds for sums of two Wigner chaoses of opposite parity under mirror symmetry. The bound contains an explicit remainder consisting of internal contractions of the two kernels. In particular, vanishing of this remainder together with the fourth free cumulant implies convergence to the standard semicircular law. If the higher-order kernel is cyclically symmetric, the internal contractions can be moved to an endpoint and we obtain an explicit fourth-root bound without a remainder. For fully symmetric kernels we obtain a square-root bound under the $W_2$-Wasserstein distance. Using the operator-norm support estimate for Wigner integrals, we also extend these bounds to every $W_p$-Wasserstein distance , where $1\le p<\infty$.
		
	\end{abstract}
	\section{Introduction}
The method of moments is a standard tool for proving Gaussian limits of nonlinear functionals of Gaussian fields \cite{GS85}. However, computing all moments can be difficult. Nualart and Peccati \cite{NP05} showed that, within a fixed Wiener chaos and under variance normalization, convergence of the fourth moment alone characterizes convergence to the standard Gaussian law. More precisely,
	let $(W_t)_{t\ge 0}$ be a standard one-dimensional Brownian motion. For any fixed integer $n\ge 1$ and deterministic square-integrable function $f$ on $\mathbb{R}^n_+$, let $I_n^W(f)$ be the $n$-th multiple Wiener integral of $f$ with respect to $W_t$; see \cite{N06}. 
	Assume $(f_k)_{k\in\mathbb{N}}$ be a sequence of real-valued symmetric functions in $L^2(\mathbb{R}_+^n)$ with $n!\left\Vert f_k\right\Vert_{L^2(\mathbb{R}_+^n)}^2=1$ for all $k$. Then the following statements are equivalent:
		\begin{enumerate}
			\item $\lim_{k\to\infty}\mathbb{E}\left[I^W_n(f_k)^4\right]=3$;
			\item The random variables $I^W_n(f_k)$ converge in distribution to the standard Gaussian distribution as $k\to\infty$.
		\end{enumerate}
Quantitative refinements using the Malliavin--Stein method are developed in \cite{NP09,NP12,NP15}. This result, known as the fourth moment theorem,  revealed a remarkable rigidity in the structure of Wiener chaoses: the Gaussian law is the only possible limit for sequences of multiple integrals with a fixed order, provided that their fourth moments converge to 3.\par
 Free Probability Theory is a branch of mathematics established by Dan Voiculescu in the 1980s. It originally emerged from research in operator algebra theory specifically within von Neumann algebras. And now free probability has evolved into an important field intersecting probability theory, random matrix theory, and combinatorial mathematics. Voiculescu's original motivation for introducing free probability was to analyze the von Neumann algebras associated with free groups. Voiculescu \cite{VOI91} proved the central limit theorem for non-commutative random variables and established that the semicircular law holds the same status in free probability theory as the normal distribution does in classical probability theory. More precisely, let $(A,\tau)$ be a non-commutative probability space. 
 Suppose $\left\{f_{n}\right\}_{n\in\mathbb{N}}$ is a sequence of freely independent self-adjoint random variables in a tracial $W^*$-probability space $(A,\tau)$ such that $\tau(f_n)=0$ for all $n$. And for any positive integer $k$, we have $\sup_{n\in\mathbb{N}}\left|\tau(f_n^k)\right|<\infty$ and
 \begin{equation*}
 	\lim_{N\to\infty}\frac{1}{N}\sum_{1\le p\le N}\tau(f_p^2)=\frac{\alpha^2}{4}>0.
 \end{equation*}
 Then
 \begin{equation*}
 	\lim_{N\to\infty}\tau(S_N^k)=\int_{-\alpha}^{\alpha}t^kd\mu_t,
 \end{equation*}
 where $\displaystyle S_N=\frac{1}{\sqrt{N}}\sum_{p=1}^{N}f_p$ is the partial summation of $\left\{f_n\right\}$ and $d\mu_t$ is probability measure on $\mathbb{R}$ with the density $\displaystyle\dfrac{2}{\pi\alpha^2}\sqrt{\alpha^2-t^2}$.Within the framework of free probability theory, Kemp, Nourdin, Peccati and Speicher \cite{KNPS12} extended the fourth-moment theorem into the case of the Wigner chaos. Let $(\mathcal{A},\tau)$ be a non-commutative probability space and $(S_t)_{t\ge 0}$ be a free Brownian motion. For any fixed integer $n\ge 1$ and deterministic mirror symmetric square-integrable function $f$ on $\mathbb{R}^n_+$, let $I_n^S(f)$ be the $n$-th multiple Wigner integral of $f$ with respect to $S_t$. 
Suppose $(f_k)_{k\in\mathbb{N}}$ is a sequence of mirror symmetric functions in $L^2(\mathbb{R}^n_+)$ with $\left\Vert f_k\right\Vert_{L^2(\mathbb{R}_+^n)}=1$ for all $k$. Then the following statements are equivalent:
		\begin{enumerate}
			\item $\lim_{k\to\infty}\tau\left[I^S_n(f_k)^4\right]=2$;
			\item The random variables $I^S_n(f_k)$ converge in distribution to the standard semicircular distribution as $k\to\infty$.
		\end{enumerate}

	Furthermore they established the convergence rate for the fourth moment theorem in the case of second order homogeneous Wigner chaos. More precisely, 
	suppose $h$ has a Fourier expansion as 
	\begin{equation*}
		h(x)=\int_{\mathbb{R}}e^{ix\xi}\nu(d\xi).
	\end{equation*}
	Define a seminorm
	\begin{equation*}
		\mathscr{S}_2(h)=\int_{\mathbb{R}}\xi^2|\nu|(d\xi).
	\end{equation*}
Let $\mathcal{C}_2$ denote the set of functions $h$ with $\mathscr{S}_2(h)<\infty$. Define
	\begin{equation*}
		d_{\mathcal{C}_2}(X,Y)=\sup\{\left|\tau\left[h(X)\right]-\tau\left[h(Y)\right]\right|:h\in\mathcal{C}_2,\mathscr{S}_2(h)\le 1\}.
	\end{equation*} 
	Then for any self-adjoint $F$ in the second homogeneous Wigner chaos with $\tau(F^2)=1$, and a standard semicircular variable $S$, one has
	\begin{equation*}
		d_{\mathcal{C}_2}(F,S)\le \frac{1}{2}\sqrt{\frac{3}{2}}\left(\tau(F^4)-2\right)^{\frac{1}{2}}.
	\end{equation*}
	\par
Bourguin and Campese \cite{BC17} extended the fourth-moment theorem of the Wigner chaos into the case of higher orders. They provided the bound 
	\begin{equation*}
		d_{\mathcal{C}_2}(F,S)\le C_n\left(\tau(F^4)-2\right)^{\frac{1}{2}}
	\end{equation*}
	for fully symmetric variables in a homogeneous Wigner chaos of order $n$ with a constant $C_n$.\par
	By the theory of free Stein discrepancy, C\'ebron \cite{CEB21} removed the assumption of fully symmetric and replaced the distance $d_{\mathcal{C}_2}$ by $W_2$ Wasserstein distance.  C\'ebron proved that for any self-adjoint element $F$ in homogeneous Wigner chaos of order $n$ with $\tau(F^2)=1$, 
	\begin{equation*}
		W_2(F,S)\le n^{\frac{3}{4}}(\tau(F^4)-2)^{\frac{1}{4}}.
	\end{equation*}
	Particularly if $F=I_n^S(f)$ for some fully symmetric $f\in L^2(\mathbb{R}_+^n)$, then
	\begin{equation*}
		W_2(F,S)\le n^{\frac{1}{2}}(\tau(F^4)-2)^{\frac{1}{2}}.
	\end{equation*}\par
	Recently Basse-O'Connor, Kramer-Bang and Svendsen \cite{CDS} first proved a fourth-moment theorem for the sum of Wiener chaoses. Suppose $q,r$ are positive integers with opposite parity. Let $f_k\in L^2(\mathbb{R}_+^q),g_k\in L^2(\mathbb{R}_+^r)$ be real-valued symmetric kernels such that $q!\left\Vert f_k\right\Vert_{L^2(\mathbb{R}_+^q)}^2+r!\left\Vert g_k\right\Vert_{L^2(\mathbb{R}_+^r)}^2=1$ for positive integer $k$. Denote $X_k=I^W_q(f_k)+I^W_r(g_k)$ the sum of Wiener chaoses. Then the following statements are equivalent
		\begin{enumerate}
			\item $\lim_{k\to\infty}\mathbb{E}\left[X_k^4\right]=3$;
			\item $X_k$ converges in distribution to the standard Gaussian distribution as $k\to\infty$.
		\end{enumerate}
		Moreover they proved that exists a constant $C_{q,r}$ such that
		\begin{equation*}
			d_{\mathrm{TV}}(X_k,\gamma)\le C_{q,r}\sqrt{\mathbb{E}\left[X_k^4\right]-3},
		\end{equation*}
		where $\gamma$ is the standard Gaussian distribution. Later,  Diez, Maia and Nourdin \cite{DMN} proved that under the same conditions $\sqrt{6}$ is an admissible constant for this total variation bound 
		with the assistance of GPT-5.\par 
For sums of two Wigner chaoses, the free Stein kernel contains mixed contractions between kernels of different orders. Kemp and Miyagawa \cite{KM25} have established a qualitative fourth-moment theorem for two chaoses of opposite parity under mirror symmetry. Here we derive a quantitative $W_2$ bound by estimating the four blocks of a free Stein kernel directly. We retain the strictly internal contractions as an explicit remainder and control the other terms by the fourth free cumulant. In particular, if this remainder tends to zero, the bound yields semicircular convergence whenever the fourth free cumulant tends to zero. For fully symmetric kernels, we obtain a sharper bound without a remainder. The qualitative theorem itself is already known from \cite{KM25}. We also give a quantitative fourth-root bound when the higher-order kernel is cyclically symmetric. \par
	For every centered self-adjoint variable $A$, we write
	\begin{equation}\label{free-fourth-cumulant}
		\kappa_4(A)=\tau\left(A^4\right)-2\tau\left(A^2\right)^2
	\end{equation}
	for its fourth free cumulant. In particular, $\kappa_4(A)=\tau(A^4)-2$ when $\tau(A^2)=1$. This convention concerns free cumulants; the corresponding classical fourth cumulant subtracts $3\mathbb{E}[A^2]^2$ instead.
	Now we state the main theorem.
	\begin{theorem}\label{Maintheorem2}
Let $S$ be a standard semicircular variable. Let $q>r\ge1$ be integers of opposite parity, and let
$f\in L^2\left(\mathbb R_+^q\right)$ and $g\in L^2\left(\mathbb R_+^r\right)$ be mirror symmetric with
$\Vert f\Vert_2^2+\Vert g\Vert_2^2=1$. Put
$X=I_q^S\left(f\right)$, $Y=I_r^S\left(g\right)$ and $F=X+Y$.
For $0\le j\le q-r$, define
\begin{equation}\label{internal-contraction-definition}
\begin{split}
D_j\left(f,g\right)\left(x_1,\ldots,x_{q-r}\right)
={}&\int_{\mathbb R_+^r}
 f\left(x_1,\ldots,x_j,t_1,\ldots,t_r,x_{j+1},\ldots,x_{q-r}\right)\\
&\quad\cdot\overline{g\left(t_1,\ldots,t_r\right)}\,dt_1\cdots dt_r,
\end{split}
\end{equation}
where an empty coordinate block is omitted. Let
\begin{equation}\label{remainder-definition}
\mathcal R_{q,r}\left(f,g\right)
=\sum_{j=1}^{q-r-1}\Vert D_j\left(f,g\right)\Vert_2^2,
\end{equation}
with the convention that this sum is zero when $q-r=1$. Then
\begin{equation}\label{main-remainder-bound}
W_2\left(F,S\right)^2
\le C_{q,r}\kappa_4\left(F\right)^{\frac{1}{2}}
 +\mathcal R_{q,r}\left(f,g\right),
\end{equation}
where
\[
A_{q,r}=\max\left\{q\sqrt{q-1},r\sqrt r+2\right\},\qquad
C_{q,r}=2\sqrt{A_{q,r}^2+q^2\left(r-1\right)}.
\]
In particular,
\[
W_2\left(F,S\right)
\le C_{q,r}^{\frac{1}{2}}\kappa_4\left(F\right)^{\frac{1}{4}}
 +\mathcal R_{q,r}\left(f,g\right)^{\frac{1}{2}}.
\]
If $f,g$ are fully symmetric, then
\[
W_2\left(F,S\right)\le\sqrt{2q}\,\kappa_4\left(F\right)^{\frac{1}{2}}.
\]

\end{theorem}

\begin{corollary}[A quantitative bound under cyclic symmetry]\label{cyclic-quantitative-corollary}
Assume the hypotheses of Theorem \ref{Maintheorem2}. Suppose, in addition, that the higher-order kernel $f$ is cyclically symmetric, that is,
\begin{equation}\label{cyclic-kernel-definition}
f\left(t_1,\ldots,t_q\right)=f\left(t_2,\ldots,t_q,t_1\right)
\end{equation}
for almost every $\left(t_1,\ldots,t_q\right)\in\mathbb R_+^q$. The lower-order kernel $g$ is required to be mirror symmetric. Define
\[
A_{q,r}^{\mathrm{cyc}}
=\max\left\{q\sqrt{q-1},r\sqrt r+q-r+1\right\},\qquad
C_{q,r}^{\mathrm{cyc}}
=2\sqrt{\left(A_{q,r}^{\mathrm{cyc}}\right)^2+q^2\left(r-1\right)}.
\]
Then
\begin{equation}\label{cyclic-quantitative-bound}
W_2\left(F,S\right)
\le\left(C_{q,r}^{\mathrm{cyc}}\right)^{\frac{1}{2}}
\kappa_4\left(F\right)^{\frac{1}{4}}.
\end{equation}
For fixed $q,r$ and a normalized sequence satisfying these symmetry assumptions,
$\kappa_4\left(F_N\right)\to0$ implies $W_p\left(F_N,S\right)\to0$ for every fixed $1\le p<\infty$.
\end{corollary}

	The paper is organized as follows. In Section 2, we introduce the basic notions of free probability and free Malliavin calculus. In Section 3, we prove Theorem \ref{Maintheorem2} and derive the corresponding $W_p$ bounds from the support of Wigner chaoses.
	\section{Notation and Basic Definition of Free Malliavin Calculus}
	\subsection{Wigner chaos}
	Let $H_{\mathbb R}=L^2_{\mathbb R}(\mathbb R_+)$, and let $H=L^2(\mathbb R_+)$ be its complexification. Let $H^{\otimes 0}=\mathbb{C}\Omega$ be a one dimensional Hilbert space generated by a unit vector $\Omega$. Let
	\begin{equation*}
		F(H)=\oplus_{n=0}^{\infty}H^{\otimes n}
	\end{equation*}
	be the Fock space associated to $H$. For each $h\in H$, let $l_h$ and $l_h^*$ be the left annihilation and creation operators defined as
	\begin{align*}
		l_h\Omega&=0,\\
		l_h(h_1\otimes\cdots\otimes h_n)&=\langle h_1,h\rangle h_2\otimes\cdots\otimes h_n,\\
		l_h^*(h_1\otimes\cdots \otimes h_n)&=h\otimes h_1\otimes\cdots\otimes h_n.
	\end{align*}
	For each $h\in L^2_{\mathbb{R}}(\mathbb{R}_+)$ let $X(h)=l_h+l_h^*$.  Let $\mathcal{S}(H)$ be the von Neumann algebra of operators on $F(H)$ generated by $\{X(h):h\in L^2_{\mathbb{R}}(\mathbb{R}_+)\}$ and let
	\begin{equation*}
		\tau(X)=\langle X\Omega,\Omega\rangle
	\end{equation*}
	for any $X\in \mathcal{S}(H)$. Then state $\tau$ is a faithful normal trace on $\mathcal{S}(H)$ and $\{X(h):h\in L^2_{\mathbb{R}}(\mathbb{R}_+)\} $ form a semi-circular system in the noncommutative probability space $(\mathcal{S}(H),\tau)$. $X(1_{[0,t]})$ is a free Brownian motion with respect to $(\mathcal{S}(H),\tau)$.
	For convenience, let $H=L^2(\mathbb{R}_+)$ and $\mathcal{A}=\mathcal{S}(H)$. 
	Let $S=(S_t)_{t\ge 0}$ be a free Brownian motion. Let $f=1_{[s_1,t_1]\times [s_2,t_2]\times\cdots\times [s_n,t_n]}$ be the off-diagonal step function which is the indicator function of $n$ dimensional rectangle such that the intervals $[s_1,t_1],[s_2,t_2],\cdots,[s_n,t_n]$ are pairwise disjoint. The Wigner integral $I_n^S(f)$ is defined by
	\begin{equation*}
		I_n^S(f)=(S_{t_1}-S_{s_1})(S_{t_2}-S_{s_2})\cdots (S_{t_n}-S_{s_n}).
	\end{equation*} 
	Extend $I_n^S$ linearly over the set of all off-diagonal step functions which is dense in $L^2(\mathbb{R}^n_+)$. 
	The freeness of the increments of $S$ yields the Wigner isometry
	\begin{equation*}
		\tau\left[I_n^S(g)^*I_n^S(f)\right]=\left\langle f,g\right\rangle_{L^2(\mathbb{R}^n_+)}.
	\end{equation*}
	which means that $I_n^S$ can be extended to the full space $L^2(\mathbb{R}_+^n)$. If $f$ is any function in $L^2(\mathbb{R}_+^n)$, then we write
	\begin{equation*}
		I_n^S(f)=\int f(t_1,\cdots,t_n)dS_{t_1}\cdots dS_{t_n}.
	\end{equation*}
	However the increments of $S$ is noncommutative. If we permute the variable of $f$, then $I_n^S(f)$ change. The image of $I_n^S$ on $L^2(\mathbb{R}_+^n)$ is called the $n$-th Wigner chaos. We use $I_0^S(c)=c1_{\mathcal A}$ and omit the superscript $S$ when no confusion is possible. By density argument we know that different orders of Wigner chaos are orthogonal and 
	\begin{equation*}
		\tau\left[I_m^S\left(g\right)^*I_n^S\left(f\right)\right]
=\begin{cases}
\left\langle f,g\right\rangle_{L^2\left(\mathbb R_+^n\right)},&m=n,\\
0,&m\ne n.
\end{cases}
	\end{equation*}
	
	Now we consider the self-adjoint elements.
	\begin{definition}
		Let $f\in L^2(\mathbb{R}_+^n)$.
		\begin{enumerate}
			\item The adjoint of $f$ is the function $\displaystyle f^*(t_1,t_2,\cdots,t_n)=\overline{f(t_n,\cdots,t_1)}$.
			\item $f$ is called mirror symmetric if $f=f^*$ for almost all $t_1,t_2,\cdots,t_n\ge 0$ with respect to the product Lebesgue measure.
			\item $f$ is called fully symmetric if it is real-valued and for any permutation $\sigma$ in the symmetric group
			\begin{equation*}
				f(t_1,t_2,\cdots,t_n)=f(t_{\sigma(1)},t_{\sigma(2)},\cdots,t_{\sigma(n)})
			\end{equation*} 
			for almost all $t_1,t_2,\cdots,t_n\ge 0$ with respect to the product Lebesgue measure.
		\end{enumerate}
	\end{definition}
	\begin{definition}\cite{BIA98}
		Let $q,r$ be two positive integers. Suppose $f\in L^2(\mathbb{R}^q_+)$ and $g\in L^2(\mathbb{R}^r_+)$, and $0\le k\le\min(q,r)$, then the $k$-th contraction of $f$ and $g$ is the function $h\in L^2(\mathbb{R}_+^{q+r-2k})$ such that
		\begin{align*}
			h&(t_1,\cdots,t_{q+r-2k})=f\stackrel{k}{\frown}g(t_1,t_2,\cdots,t_{q+r-2k})\\
			&=\int_{\mathbb{R}_+^k}f(t_1,\cdots,t_{q-k},s_k,\cdots,s_1)g(s_1,\cdots,s_k,t_{q-k+1},\cdots,t_{q+r-2k})ds_1\cdots ds_k.
		\end{align*}
	\end{definition}
	In \cite{BIA98}, Biane and Speicher established the following contraction formula for the Wigner chaos. 
	\begin{theorem}\label{contraction formula 1}\cite{BIA98}
		Given $f\in L^2(\mathbb{R}_+^q),g\in L^2(\mathbb{R}_+^r)$, then
		\begin{equation*}
			I_q(f)I_r(g)=\sum_{k=0}^{\min(q,r)}I_{q+r-2k}(f\stackrel{k}{\frown}g)
		\end{equation*}
	\end{theorem}
	Fix a positive integer $n$ and positive integers $n_1,\ldots,n_k$ satisfying $n_1+\cdots+n_k=n$. The corresponding consecutive subsets of $\{1,\ldots,n\}$ form the block decomposition denoted by $(n_1,\ldots,n_k)$. Let $\mathcal{P}_2(n_1\otimes\cdots\otimes n_k)$ be the set of pairings $\pi$ of $\{1,\ldots,n\}$ such that no pair contains more than one element from any given block. A pairing is called noncrossing if there is no quadruple $(i,j,k,l)$ with $1\le i<j<k<l\le n$ such that $\{i,k\}\in\pi$ and $\{j,l\}\in\pi$. The set of all such noncrossing pairings is denoted by $NC_2(n_1\otimes\cdots\otimes n_k)$. These pairings are frequently used to compute traces of cross terms.
	The following identity provides a computational method for the trace of products of several Wigner chaoses.
	\begin{theorem}\cite{KNPS12}\label{product formula for Wigner integral 2}
		Let $n_1,n_2,\cdots,n_l$ be positive integers, and suppose $f_1,\cdots,f_l$ are functions with $f_i\in L^2(\mathbb{R}_+^{n_i})$ for $1\le i\le l$. Then 
		\begin{equation*}
			\tau\left(I_{n_1}(f_1)\cdots I_{n_l}(f_l)\right)=\sum_{\pi\in NC_2(n_1\otimes \cdots\otimes n_l)}\int_{\pi}f_1\otimes\cdots\otimes f_l.
		\end{equation*}
		Here $NC_2(n_1\otimes \cdots\otimes n_l)$ is the set of all noncrossing pairing of $(n_1,\cdots,n_l)$.
	\end{theorem}
	Particularly the following formula establishes the connection between the fourth moment of a Wigner chaos and its contractions.
	\begin{theorem}\cite{KNPS12}\label{integration esitimate}
		Let $f=f^*\in L^2(\mathbb{R}_+^n)$ such that $\left\Vert f\right\Vert_{L^2(\mathbb{R}_+^n)}=1$. We have
		\begin{equation*}
			\tau\left[\left(I_n(f)\right)^4\right]=2+\sum_{k=1}^{n-1}\left\Vert f\stackrel{k}{\frown} f\right\Vert^2_{L^2(\mathbb{R}_+^{2n-2k})}.
		\end{equation*} 
	\end{theorem}
	In order to estimate the integral with respect to pairing, we need the following lemma.
	\begin{lemma}\cite{KNPS12}\label{estimate of pairing integration}
		Let $n_1,n_2,\cdots,n_k$ be positive integer and $f_i\in L^2(\mathbb{R}_+^{n_i})$, then for any pairing $\pi\in\mathcal{P}_2(n_1\otimes n_2\otimes\cdots\otimes n_k)$ we have
		\begin{equation*}
			\left|\int_{\pi}f_1\otimes \cdots\otimes f_k\right|\le \prod_{i=1}^{k}\left\Vert f_i\right\Vert_{L^2(\mathbb{R}_+^{n_i})}.
		\end{equation*}
	\end{lemma}
	Due to noncommutativity, it is also necessary to introduce the Wigner chaos in the form of a second-order tensor to give a representation of square-integrable biprocesses. 
	For nonnegative integers $m,n$ and $f\in  L^2(\mathbb{R}_+^m),g\in L^2(\mathbb{R}_+^n)$. Let $h=f\otimes g$ and define
	\begin{equation*}
		\left[I_m^S\otimes I_n^S\right](h)=\int h(t_1,\cdots,t_m;s_1,\cdots,s_n)dS_{t_1}\cdots dS_{t_m}\otimes dS_{s_1}\cdots dS_{s_n}=I_m^S(f)\otimes I_n^S(g).
	\end{equation*}
	And in general $I_m^S\otimes I_n^S$ extend to all function in $L^2(\mathbb{R}_+^m)\otimes L^2(\mathbb{R}_+^n)=L^2(\mathbb{R}_+^{m+n})$,
	\begin{equation*}
		I_m^S\otimes I_n^S: L^2(\mathbb{R}_+^{m+n})\to L^2(\mathcal A\otimes\mathcal A,\tau\otimes\tau).
	\end{equation*}
	Write $I_{m,n}^S=[I_m^S\otimes I_n^S]$. By a density argument,
	\begin{equation*}
		\tau\otimes\tau\left[I_{m,n}^S\left(g\right)^*I_{q,r}^S\left(f\right)\right]
=\begin{cases}
\left\langle f,g\right\rangle_{L^2\left(\mathbb R_+^{m+n}\right)},&\left(m,n\right)=\left(q,r\right),\\
0,&\left(m,n\right)\ne\left(q,r\right).
\end{cases}
	\end{equation*}
	Now for any bi-process $I_m^S\otimes I_n^S(f_t)$ define
	\begin{equation*}
		\int I_m^S\otimes I_n^S(f_t)\# dS_t=\int f_t(t_1,\cdots,t_m,s_1,\cdots,s_n)dS_{t_1}dS_{t_2}\cdots dS_{t_m}dS_tdS_{s_1}\cdots dS_{s_n}.
	\end{equation*}
	In \cite{CEB21}, C\'ebron proved the following contraction formula for the Wigner bi-chaos.
	\begin{theorem}\cite{CEB21}\label{product formula for bi-integral 1}
		Let $1\le k\le n$. For all $f,g\in L^2(\mathbb{R}_+^{n-1})$, we have
		\begin{equation*}
			I_{n-1}(f)\cdot \left[ I_{k-1}\otimes I_{n-k}\right](g)=\sum_{i=0}^{k-1}\left[I_{n-1+k-1-2i}\otimes I_{n-k}\right](f\stackrel{i}{\frown}g).
		\end{equation*}
	\end{theorem}
More generally, if $f\in L^2\left(\mathbb R_+^{q-1}\right)$,
$g\in L^2\left(\mathbb R_+^{r-1}\right)$ and $1\le k\le r$, then
\begin{equation}\label{mixed-bi-integral-product}
I_{q-1}\left(f\right)\cdot\left[I_{k-1}\otimes I_{r-k}\right]\left(g\right)=\sum_{i=0}^{\min\left(q-1,k-1\right)}
 \left[I_{q+k-2i-2}\otimes I_{r-k}\right]\left(f\stackrel i\frown g\right).
\end{equation}
Here multiplication by $I_{q-1}\left(f\right)$ acts on the first tensor leg.

	\subsection{Free Malliavin calculus and free Stein discrepancy}
	In \cite{BIA98}, Biane and Speicher defined the free Malliavin derivative. Let $\mathcal S_{alg}(H)$ denote the algebra of noncommutative polynomials generated by $X(h)$, $h\in H_{\mathbb R}$.
	\begin{definition}
		The mapping
		\begin{equation*}
			\nabla:\mathcal S_{alg}(H)\to\mathcal S_{alg}(H)\otimes H\otimes\mathcal S_{alg}(H),
		\end{equation*}
		with the properties 
		\begin{enumerate}
			\item For all $A,B\in\mathcal S_{alg}(H)$, 
			\begin{equation*}
				\nabla(AB)=A\cdot(\nabla B)+(\nabla A)\cdot B.
			\end{equation*}
			\item For all $h_1,h_2,\cdots,h_n\in H_{\mathbb R}$,
			\begin{equation*}
				\nabla(X(h_1)X(h_2)\cdots X(h_n))=\sum_{k=1}^{n}(X(h_1)\cdots X(h_{k-1}))\otimes h_k\otimes (X(h_{k+1})\cdots X(h_n)).
			\end{equation*}
		\end{enumerate}
		is called the free Malliavin derivative.
	\end{definition}
The derivation is extended complex linearly and then closed in the corresponding spaces. We use the same notation for its closure, whose domain contains each finite sum of Wigner chaoses.
	More precisely when the free Malliavin derivative acts on Wigner integral we have the following theorem.
	\begin{theorem}\label{Malliavin derivative of Wigner chaos}\cite{BIA98}
		Let $n\ge1$, let $\mathcal{H}_r$ denote the $r$-th homogeneous Wigner chaos, and put $\mathcal{H}_0=\mathbb{C}1_{\mathcal{A}}$. The Malliavin derivative $\nabla$ maps $\mathcal{H}_n$ into $L^2\left(\mathbb{R}_+;\bigoplus_{a+b=n-1}\mathcal{H}_a\otimes\mathcal{H}_b\right)$. More precisely, if $f\in L^2(\mathbb{R}_+^n)$, then for almost all $t\ge 0$
		\begin{equation*}
			\nabla_t\left(I_n\left(f\right)\right)=\sum_{k=1}^{n}\left(I_{k-1}\otimes I_{n-k}\right)\left(f_t^k\right),
		\end{equation*}
		where $f_t^k(t_1,\cdots,t_{k-1},t_{k+1},\cdots,t_n)=f(t_1,\cdots,t_{k-1},t,t_{k+1},\cdots,t_n)$.
	\end{theorem}
	In order to give an improvement of the logarithmic Sobolev inequality of free diffusion. Fathi and Nelson \cite{FB17} introduced the free Stein discrepancy.
	Let $\partial$ be the free difference quotient on polynomials, defined by
\[
\partial 1=0,\qquad
\partial z^m=\sum_{j=0}^{m-1}z^j\otimes z^{m-1-j},\qquad m\ge1,
\]
and extended complex linearly. Here $\partial P\left(F\right)$ denotes the evaluation of both tensor factors of $\partial P$ at $F$.
	\begin{definition}\cite{FB17}
		Let $F$ be a self-adjoint element in $(A,\tau)$. A free Stein kernel of $F$ is a element $\xi \in L^2(A,\tau)\otimes L^2(A,\tau)$ such that for any $P\in \mathbb{C}\langle X\rangle$
		\begin{equation*}
			\left\langle F,P(F)\right\rangle_{L^2(\mathcal{A},\tau)}=\left\langle\xi ,\partial P(F)\right\rangle_{L^2(A,\tau)\otimes L^2(A,\tau)}.
		\end{equation*}
		The free Stein discrepancy of $F$ with respect to semicircular varible $S$ is defined by
		\begin{equation*}
			\Sigma^*(F|S)=\inf_{\xi}\Vert \xi-1_{\mathcal{A}}\otimes 1_{\mathcal{A}}\Vert_{L^2(A,\tau)\otimes L^2(A,\tau)},
		\end{equation*} 
		where the infimum is taken over all free Stein kernels of $F$.
	\end{definition}
	In \cite{CEB21}, C\'ebron proved the following free WS inequality.
	\begin{proposition}\cite{CEB21}
		For all self-adjoint element $F$ in $(\mathcal{A},\tau)$, we have
		\begin{equation*}
			W_2(F,S)\le \Sigma^*(F|S).
		\end{equation*} 
	\end{proposition}
	The following construction applies to finite sums of Wigner chaoses.
	\begin{proposition}\cite[Corollary 3.10 and Proposition 3.11]{CEB21}\label{finite-chaos-Stein-kernel}
		Let $F$ be centered and self-adjoint with a finite Wigner chaos expansion. Then
		\begin{equation*}
			K_F=\int_0^\infty\left[\left(\operatorname{id}\otimes\tau\right)\left(\nabla_tF\right)\otimes1_{\mathcal{A}}\right]\left(\nabla_tF\right)^*\,dt
		\end{equation*}
		belongs to $L^2(\mathcal{A},\tau)\otimes L^2(\mathcal{A},\tau)$ and is a free Stein kernel of $F$. Consequently,
		\begin{equation*}
			W_2(F,S)\le\Sigma^*(F|S)\le\left\Vert K_F-1_{\mathcal{A}}\otimes1_{\mathcal{A}}\right\Vert_2.
		\end{equation*}
		If $F=I_n^S(f)$, $n\ge2$, and $\tau(F^2)=1$, then
		\begin{equation*}
			\left\Vert K_F-1_{\mathcal{A}}\otimes1_{\mathcal{A}}\right\Vert_2\le n^{\frac{3}{4}}\kappa_4(F)^{\frac{1}{4}}.
		\end{equation*}
		Here $S$ is standard semicircular, and $\Vert\cdot\Vert_2$ is the tensor-product Hilbert norm. In the sequel, multiplication by $\left(\operatorname{id}\otimes\tau\right)(\nabla_tF)$ means multiplication on the first tensor leg.
	\end{proposition}

\section{Proof of Theorem \ref{Maintheorem2}}

\begin{proof}[Proof of Theorem \ref{Maintheorem2} and Corollary \ref{cyclic-quantitative-corollary}]
All the $L^2$-norms below are taken in the corresponding kernel spaces, or in the tensor-product Hilbert space when applied to biprocesses. Since $q\ne r$, orthogonality and the Wigner isometry give
\[
\tau\left(F\right)=0,\qquad
\tau\left(F^2\right)=\tau\left(X^2\right)+\tau\left(Y^2\right)
=\Vert f\Vert_2^2+\Vert g\Vert_2^2=1.
\]
The kernel calculations may first be carried out for finite linear combinations of elementary tensors and then extended by $L^2$ approximation and continuity of contractions. We first compute the fourth moment of $F$. Since $q,r$ have opposite parity, the total orders in each of $X^3Y$, $XY^3$, $XY$ and $YX$ are odd. There is therefore no corresponding noncrossing pairing, and Theorem \ref{product formula for Wigner integral 2} yields
\[
\tau\left(X^3Y\right)=\tau\left(XY^3\right)
=\tau\left(XY\right)=\tau\left(YX\right)=0.
\]
Expanding the fourth power and using traciality, we obtain
\begin{align*}
\tau\left[\left(X+Y\right)^4\right]
&=\tau\left(X^4\right)+4\tau\left(X^3Y\right)+4\tau\left(X^2Y^2\right)
 +2\tau\left(XYXY\right)+4\tau\left(XY^3\right)+\tau\left(Y^4\right)\\
&=\tau\left(X^4\right)+4\tau\left(X^2Y^2\right)
 +2\tau\left(XYXY\right)+\tau\left(Y^4\right).
\end{align*}
By the contraction formula,
\[
XY=I_q\left(f\right)I_r\left(g\right)
=\sum_{k=0}^r I_{q+r-2k}\left(f\stackrel{k}{\frown}g\right).
\]
Combining this identity with orthogonality gives
\begin{align*}
\tau\left(XYXY\right)
&=\sum_{k=0}^r\tau\left[\left(I_{q+r-2k}\left(f\stackrel{k}{\frown}g\right)\right)^2\right]\\
&=\operatorname{Re}\sum_{k=0}^r
 \left\langle f\stackrel{k}{\frown}g,g\stackrel{k}{\frown}f\right\rangle_{L^2\left(\mathbb R_+^{q+r-2k}\right)}.
\end{align*}
Similarly,
\begin{align*}
\tau\left(X^2Y^2\right)
&=\sum_{i=0}^q\sum_{j=0}^r
 \tau\left[I_{2q-2i}\left(f\stackrel{i}{\frown}f\right)
 I_{2r-2j}\left(g\stackrel{j}{\frown}g\right)\right]\\
&=\sum_{i=0}^r
 \tau\left[I_{2r-2i}\left(f\stackrel{i+q-r}{\frown}f\right)
 I_{2r-2i}\left(g\stackrel{i}{\frown}g\right)\right].
\end{align*}
Also, self-adjointness of $X,Y$, the trace property and the Wigner isometry imply
\begin{equation}\label{covariance-id}
\tau\left(X^2Y^2\right)=\tau\left(YX^2Y\right)
=\sum_{i=0}^r\left\Vert f\stackrel{i}{\frown}g\right\Vert_{L^2\left(\mathbb R_+^{q+r-2i}\right)}^2.
\end{equation}
The fourth-moment formula for a single Wigner integral gives
\begin{equation}\label{individual-cumulants}
\begin{aligned}
\kappa_4\left(X\right)&=\sum_{l=1}^{q-1}\left\Vert f\stackrel l\frown f\right\Vert_2^2,\\
\kappa_4\left(Y\right)&=\sum_{l=1}^{r-1}\left\Vert g\stackrel l\frown g\right\Vert_2^2.
\end{aligned}
\end{equation}
By mirror symmetry and reversal of the integration variables,
\[
D_{q-r}\left(f,g\right)=f\stackrel r\frown g,\qquad
D_0\left(f,g\right)=\left(f\stackrel r\frown g\right)^*.
\]
To compute the inner product of $f\otimes g$ with its adjoint, let $x,s$ each denote $r$ coordinates and let $t$ denote $q-r$ coordinates. Fubini's theorem gives
\begin{align*}
\left\langle f\otimes g,g\otimes f\right\rangle_2
&=\int f\left(x,t\right)g\left(s\right)\overline{g\left(x\right)}\,
 \overline{f\left(t,s\right)}\,dx\,dt\,ds\\
&=\int D_0\left(f,g\right)\left(t\right)
 \overline{\left(f\stackrel r\frown g\right)\left(t\right)}\,dt\\
&=\overline{\left\langle f\stackrel r\frown g,
 \left(f\stackrel r\frown g\right)^*\right\rangle_2}.
\end{align*}
In particular, the two inner products have the same real part. Combining this identity with the fourth-moment expansion above, and using $\tau\left(F^2\right)=1$ and
$\Vert f\otimes g\Vert_2^2=\Vert f\Vert_2^2\Vert g\Vert_2^2$, we obtain
\begin{equation}\label{positive-mixed-decomposition}
\begin{aligned}
&\kappa_4\left(F\right)-\kappa_4\left(X\right)-\kappa_4\left(Y\right)\\
&=4\sum_{l=1}^r\left\Vert f\stackrel l\frown g\right\Vert_2^2
 +2\operatorname{Re}\sum_{l=0}^r
 \left\langle f\stackrel l\frown g,
 \left(f\stackrel l\frown g\right)^*\right\rangle_2\\
&=2\sum_{l=1}^{r-1}\left\Vert f\stackrel l\frown g\right\Vert_2^2
 +\sum_{l=1}^r\left\Vert f\stackrel l\frown g+
 \left(f\stackrel l\frown g\right)^*\right\Vert_2^2\\
&\quad+\left\Vert f\stackrel r\frown g+
 \left(f\stackrel r\frown g\right)^*\right\Vert_2^2\ge0.
\end{aligned}
\end{equation}
Consequently,
\begin{equation}\label{cumulants-below-total-fourth-moment}
\kappa_4\left(X\right)+\kappa_4\left(Y\right)\le\kappa_4\left(F\right).
\end{equation}

We next compute the free Stein kernel. By the formula for the free Malliavin derivative in Theorem \ref{Malliavin derivative of Wigner chaos},
\[
\nabla_tF=\sum_{k=1}^q\left(I_{k-1}\otimes I_{q-k}\right)\left(f_t^k\right)
 +\sum_{k=1}^r\left(I_{k-1}\otimes I_{r-k}\right)\left(g_t^k\right).
\]
Taking the partial trace leaves only the terms whose right tensor leg has order zero, so
\[
\left(\operatorname{id}\otimes\tau\right)\left(\nabla_tF\right)
=I_{q-1}\left(f_t^q\right)+I_{r-1}\left(g_t^r\right).
\]
Since $\left(A\otimes B\right)^*=A^*\otimes B^*$, taking adjoints in a Wigner bi-integral gives
	\begin{align*}
		&\left(\int_{\mathbb{R}_+^{n-1}}f\left(t_1,\cdots,t_{k-1},t_{k+1},\cdots,t_n\right)dS_{t_1}\cdots dS_{t_{k-1}}\otimes dS_{t_{k+1}}\cdots dS_{t_n}\right)^*\\
		&=\int_{\mathbb{R}_+^{n-1}}\overline{f\left(t_{k-1},\cdots,t_1,t_n,\cdots,t_{k+1}\right)}dS_{t_1}\cdots dS_{t_{k-1}}\otimes dS_{t_{k+1}}\cdots dS_{t_n}.
	\end{align*}
Write
	\begin{equation*}
		\widetilde{f}_t^k\left(t_1,\cdots,t_{k-1},t_{k+1},\cdots,t_n\right)=\overline{f\left(t_{k-1},\cdots,t_1,t,t_n,\cdots,t_{k+1}\right)},
	\end{equation*}
The same notation is used for $g$. When all coordinates are displayed, $\widetilde f$ and $\widetilde g$ denote these blockwise reversed and conjugated kernels at the fixed tensor cut under consideration. The adjoint of a Wigner bi-integral can then be expressed as
	\begin{equation*}
		\left(\left[I_{k-1}\otimes I_{n-k}\right]\left(f_t^k\right)\right)^*=\left[I_{k-1}\otimes I_{n-k}\right]\left(\widetilde{f}_t^k\right).
	\end{equation*}
Applying the contraction formula from Theorem \ref{contraction formula 1} on the first tensor leg, we obtain
	\begin{align*}
		&\left(\operatorname{id}\otimes\tau\right)\left(\nabla_tF\right)\cdot \left(\nabla_tF\right)^*\\
		&=\left(I_{q-1}\left(f^{q}_t\right)+I_{r-1}\left(g^{r}_t\right)\right)\cdot \left(\sum_{k=1}^{q}\left(I_{k-1}\otimes I_{q-k}\right)\left(\widetilde{f}_t^k\right)+\sum_{k=1}^{r}\left(I_{k-1}\otimes I_{r-k}\right)\left(\widetilde{g}_t^k\right)\right)\\
		&=\sum_{k=1}^{q}\sum_{i=0}^{k-1}\left[I_{q-1+k-1-2i}\otimes I_{q-k}\right]\left(f_t^q\stackrel{i}{\frown}\widetilde{f}_t^k\right)\\
&\quad+\sum_{k=1}^{q}\sum_{i=0}^{\min\left(k-1,r-1\right)}\left[I_{r-1+k-1-2i}\otimes I_{q-k}\right]\left(g_t^r\stackrel{i}{\frown}\widetilde{f}_t^k\right)\\
		&\quad +\sum_{k=1}^{r}\sum_{i=0}^{k-1}\left[I_{q-1+k-1-2i}\otimes I_{r-k}\right]\left(f_t^q\stackrel{i}{\frown}\widetilde{g}_t^k\right)\\
&\quad+\sum_{k=1}^{r}\sum_{i=0}^{k-1}\left[I_{r-1+k-1-2i}\otimes I_{r-k}\right]\left(g_t^r\stackrel{i}{\frown}\widetilde{g}_t^k\right).
	\end{align*}
We now identify the terms in the four sums with the same ordered pair of chaos degrees. There are two types of overlap:
	\begin{enumerate}
\item If $q-r=k_1-k_2$ and $i_1=i_2$, then the terms
		\begin{equation*}
			\left[I_{r-1+k_1-1-2i_1}\otimes I_{q-k_1}\right]\left(g_t^r\stackrel{i_1}{\frown}\widetilde{f}_t^{k_1}\right),\ \left[I_{q-1+k_2-1-2i_2}\otimes I_{r-k_2}\right]\left(f_t^q\stackrel{i_2}{\frown}\widetilde{g}_t^{k_2}\right)
\end{equation*} have the same ordered pair of chaos degrees;
\item If $q-r=i_1-i_2=k_1-k_2$, then the terms
		\begin{equation*}
\begin{gathered}
\left[I_{r+k-2i-2}\otimes I_{r-k}\right]\left(f_t^q\stackrel{q-r+i}\frown\widetilde f_t^{q-r+k}\right),\\
\left[I_{r+k-2i-2}\otimes I_{r-k}\right]\left(g_t^r\stackrel i\frown\widetilde g_t^k\right)
\end{gathered}
\end{equation*} have the same ordered pair of chaos degrees.
	\end{enumerate}
Consider the first case. Here $0\le i\le r-1$ and $q-r+i+1\le k\le q$. Combining the kernels gives
	\begin{equation*}
		\left[I_{r-1+k-1-2i}\otimes I_{q-k}\right]\left(g_t^r\stackrel{i}{\frown}\widetilde{f}_t^{k}+f_t^q\stackrel{i}{\frown}\widetilde{g}_t^{k-q+r}\right).
	\end{equation*}
The corresponding squared $L^2$-norm expands as
	\begin{align}
		&\left\Vert\int_{\mathbb{R}_+}\left(g_t^r\stackrel{i}{\frown}\widetilde{f}_t^{k}+f_t^q\stackrel{i}{\frown}\widetilde{g}_t^{k-q+r}\right)dt\right\Vert_{L^2\left(\mathbb{R}_+^{q+r-2-2i}\right)}^2\notag\\
		& =	\left\Vert\int_{\mathbb{R}_+}\left(g_t^r\stackrel{i}{\frown}\widetilde{f}_t^{k}\right)dt\right\Vert_{L^2\left(\mathbb{R}_+^{q+r-2-2i}\right)}^2+	\left\Vert\int_{\mathbb{R}_+}\left(f_t^q\stackrel{i}{\frown}\widetilde{g}_t^{k-q+r}\right)dt\right\Vert_{L^2\left(\mathbb{R}_+^{q+r-2-2i}\right)}^2\notag\\
		&\quad +2\operatorname{Re}\left\langle\int_{\mathbb{R}_+}\left(g_t^r\stackrel{i}{\frown}\widetilde{f}_t^{k}\right)dt,\int_{\mathbb{R}_+}\left(f_t^q\stackrel{i}{\frown}\widetilde{g}_t^{k-q+r}\right)dt\right\rangle_{L^2\left(\mathbb{R}_+^{q+r-2-2i}\right)}.
	\end{align}
By the definition of the contraction $\stackrel{i}{\frown}$,
	\begin{align*}
		&g_t^r\stackrel{i}{\frown}\widetilde{f}_t^{k}\\
		&=\int_{\mathbb{R}_+^i} g\left(x_1,\cdots,x_{r-i-1},u_i,\cdots,u_1,t\right)\\
		&\quad \widetilde{f}\left(u_1,\cdots,u_i,x_{r-i},\cdots,x_{r+k-2i-2},t,x_{r+k-2i-1},\cdots,x_{q+r-2-2i}\right)du_1\cdots du_i\\
		&=\int_{\mathbb{R}_+^i} g\left(x_1,\cdots,x_{r-i-1},u_i,\cdots,u_1,t\right)\\&
\quad \overline{f\left(x_{r+k-2i-2},\cdots,x_{r-i},u_i,\cdots,u_1,t,x_{q+r-2-2i},\cdots,x_{r+k-2i-1}\right)}du_1\cdots du_i,\\
	\end{align*}
We next compute the squared $L^2$-norm of $\int_{\mathbb{R}^+}g_t^r\stackrel{i}{\frown}\widetilde f_t^k$:
	\begin{align*}
		&\left\Vert\int_{\mathbb{R}_+}\left(g_t^r\stackrel{i}{\frown}\widetilde{f}_t^{k}\right)dt\right\Vert_{L^2\left(\mathbb{R}_+^{q+r-2-2i}\right)}^2\\
		&=\int_{\mathbb R_+^{q+r}}g\left(x_1,\cdots,x_{r-i-1},u_i,\cdots,u_1,t\right)\\
		&\quad \overline{f\left(x_{r+k-2i-2},\cdots,x_{r-i},u_i,\cdots,u_1,t,x_{q+r-2-2i},\cdots,x_{r+k-2i-1}\right)} \\
		&\quad \overline{g\left(x_1,\cdots,x_{r-i-1},u_i',\cdots,u_1',t'\right)}\\
		&\quad f\left(x_{r+k-2i-2},\cdots,x_{r-i},u'_i,\cdots,u'_1,t',x_{q+r-2-2i},\cdots,x_{r+k-2i-1}\right)\\
		&\quad du_1\cdots du_idtdu_1'\cdots du_i'dt'dx_1\cdots dx_{q+r-2-2i}\\
		&=\int_{\mathbb R_+^{q+r}}g\left(x_1,\cdots,x_{r-i-1},u_i,\cdots,u_1,t\right)\\
		&\quad f\left(x_{r+k-2i-1},\cdots,x_{q+r-2-2i},t,u_1,\cdots,u_i,x_{r-i},\cdots,x_{r+k-2i-2}\right) \\
		&\quad g\left(t',u_1',\cdots,u_{i}',x_{r-i-1},\cdots,x_1\right)\\
		&\quad f\left(x_{r+k-2i-2},\cdots,x_{r-i},u'_i,\cdots,u'_1,t',x_{q+r-2-2i},\cdots,x_{r+k-2i-1}\right)\\
		&\quad du_1\cdots du_idtdu_1'\cdots du_i'dt'dx_1\cdots dx_{q+r-2-2i}.\\
	\end{align*}
Integrating first with respect to $x_1,\ldots,x_{r-i-1}$ produces a pairing $\pi\in\mathcal P_2\left(q\otimes\left(2i+2\right)\otimes q\right)$ such that
	\begin{equation}
		\left\Vert\int_0^\infty\left(g_t^r\stackrel{i}{\frown}\widetilde{f}_t^{k}\right)dt\right\Vert_{L^2\left(\mathbb{R}_+^{q+r-2-2i}\right)}^2=\int_{\pi}f\otimes \left(g\stackrel{r-i-1}{\frown}g\right)\otimes f.
	\end{equation}
When $k=q$, the integral admits a further contraction. In this case,
	\begin{align*}
		&\left\Vert\int_0^\infty\left(g_t^r\stackrel{i}{\frown}\widetilde{f}_t^{q}\right)dt\right\Vert_{L^2\left(\mathbb{R}_+^{q+r-2-2i}\right)}^2\\
		&=\int_{\mathbb R_+^{q+r}}g\left(x_1,\cdots,x_{r-i-1},u_i,\cdots,u_1,t\right)\\
		&\quad f\left(t,u_1,\cdots,u_i,x_{r-i},\cdots,x_{q+r-2i-2}\right) \\
		&\quad g\left(t',u_1',\cdots,u_{i}',x_{r-i-1},\cdots,x_1\right)\\
		&\quad f\left(x_{q+r-2i-2},\cdots,x_{r-i},u_i',\cdots,u_1',t'\right)\\
		&\quad du_1\cdots du_idtdu_1'\cdots du_i'dt'dx_1\cdots dx_{q+r-2-2i}\\
	\end{align*}
Integrating with respect to $x_1,\ldots,x_{r-i-1},x_{r-i},\ldots,x_{q+r-2i-2}$ gives a pairing $\pi\in NC_2\left(\left(2i+2\right)\otimes\left(2i+2\right)\right)$ for which
	\begin{equation}
		\left\Vert\int_0^\infty\left(g_t^r\stackrel{i}{\frown}\widetilde{f}_t^{q}\right)dt\right\Vert_{L^2\left(\mathbb{R}_+^{q+r-2-2i}\right)}^2=\int_{\pi}\left(f\stackrel{q-i-1}{\frown} f\right)\otimes \left(g\stackrel{r-i-1}{\frown}g\right)
	\end{equation}
The set $NC_2\left(\left(2i+2\right)\otimes\left(2i+2\right)\right)$ contains exactly one pairing. Consequently,
	\begin{equation*}
		\int_{\pi}\left(f\stackrel{q-i-1}{\frown} f\right)\otimes \left(g\stackrel{r-i-1}{\frown}g\right)=\tau\left[I_{2i+2}\left(f\stackrel{q-i-1}{\frown}f\right)I_{2i+2}\left(g\stackrel{r-i-1}{\frown}g\right)\right].
	\end{equation*}
If $f$ and $g$ are fully symmetric, the same identity holds for every admissible $k$:
	\begin{align*}
		&\left\Vert\int_{\mathbb{R}_+}\left(g_t^r\stackrel{i}{\frown}\widetilde{f}_t^{k}\right)dt\right\Vert_{L^2\left(\mathbb{R}_+^{q+r-2-2i}\right)}^2\\
		&=\left\langle f\stackrel{q-i-1}{\frown}f,g\stackrel{r-i-1}{\frown}g\right\rangle_{L^2\left(\mathbb{R}_+^{2i+2}\right)}=\tau\left[I_{2i+2}\left(f\stackrel{q-i-1}{\frown}f\right)I_{2i+2}\left(g\stackrel{r-i-1}{\frown}g\right)\right].
	\end{align*}
The same calculation for the other mixed term gives a pairing $\pi\in\mathcal P_2\left(r\otimes\left(2i+2\right)\otimes r\right)$ satisfying
	\begin{equation}
		\left\Vert\int_0^\infty\left(f_t^q\stackrel{i}{\frown}\widetilde{g}_t^{k-q+r}\right)dt\right\Vert_{L^2\left(\mathbb{R}_+^{q+r-2-2i}\right)}^2=\int_{\pi}g\otimes \left(f\stackrel{q-i-1}{\frown}f\right)\otimes g.
	\end{equation}
In particular, when $k=q$, there is a pairing $\pi\in NC_2\left(\left(2i+2\right)\otimes\left(2i+2\right)\right)$ such that
	\begin{equation}
		\left\Vert\int_0^\infty\left(f_t^q\stackrel{i}{\frown}\widetilde{g}_t^{r}\right)dt\right\Vert_{L^2\left(\mathbb{R}_+^{q+r-2-2i}\right)}^2=\int_{\pi}\left(f\stackrel{q-i-1}{\frown}f\right)\otimes \left(g\stackrel{r-i-1}{\frown} g\right).
	\end{equation}
If $f$ and $g$ are fully symmetric, then for every admissible $k$,
	\begin{equation}
		\left\Vert\int_{\mathbb{R}_+}\left(f_t^q\stackrel{i}{\frown}\widetilde{g}_t^{k-q+r}\right)dt\right\Vert_{L^2\left(\mathbb{R}_+^{q+r-2i-2}\right)}^2=\tau\left[I_{2i+2}\left(f\stackrel{q-i-1}{\frown}f\right)I_{2i+2}\left(g\stackrel{r-i-1}{\frown}g\right)\right].
	\end{equation}
For the cross term, the full coordinate calculation is
	\begin{align*}
		&\quad \left\langle\int_{\mathbb{R}_+}\left(g_t^r\stackrel{i}{\frown}\widetilde{f}_t^{k}\right)dt,\int_{\mathbb{R}_+}\left(f_t^q\stackrel{i}{\frown}\widetilde{g}_t^{k-q+r}\right)dt\right\rangle_{L^2\left(\mathbb{R}_+^{q+r-2-2i}\right)}\\
		&=\int_{\mathbb R_+^{q+r}}g\left(x_1,\cdots,x_{r-i-1},u_i,\cdots,u_1,t\right)\\
		&\quad \overline{f\left(x_{r+k-2i-2},\cdots,x_{r-i},u_i,\cdots,u_1,t,x_{q+r-2-2i},\cdots,x_{r+k-2i-1}\right)} \\
		&\quad \overline{f\left(x_1,\cdots,x_{q-i-1},u_i',\cdots,u_1',t'\right)}\\
		&\quad \overline{\widetilde{g}\left(u_1',\cdots,u_i',x_{q-i},\cdots,x_{r+k-2i-2},t',x_{r+k-2i-1},\cdots,x_{q+r-2i-2}\right)}\\
		&\quad du_1\cdots du_idtdu_1'\cdots du_i'dt' dx_1\cdots dx_{q+r-2-2i}\\
		&=\int_{\mathbb R_+^{q+r}}g\left(x_1,\cdots,x_{r-i-1},u_i,\cdots,u_1,t\right)\\
		&\quad f\left(x_{r+k-2i-1},\cdots,x_{q+r-2-2i},t,u_1,\cdots,u_i,x_{r-i},\cdots,x_{r+k-2i-2}\right)\\
		&\quad f\left(t',u_1',\cdots,u_i',x_{q-i-1},\cdots,x_1\right)\\
		&\quad g\left(x_{r+k-2i-2},\cdots,x_{q-i},u_i',\cdots,u_1',t',x_{q+r-2-2i},\cdots,x_{r+k-2i-1}\right)\\
		&\quad du_1\cdots du_idtdu_1'\cdots du_i'dt' dx_1\cdots dx_{q+r-2-2i}. 
	\end{align*}
Integrating with respect to $x_1,\ldots,x_{r-i-1}$ yields a pairing $\pi\in\mathcal P_2\left(q\otimes\left(q-r+2i+2\right)\otimes r\right)$ for which
	\begin{align}
		&\left\langle\int_{\mathbb{R}_+}\left(g_t^r\stackrel{i}{\frown}\widetilde{f}_t^{k}\right)dt,\int_{\mathbb{R}_+}\left(f_t^q\stackrel{i}{\frown}\widetilde{g}_t^{k-q+r}\right)dt\right\rangle_{L^2\left(\mathbb{R}_+^{q+r-2-2i}\right)}\notag\\
		&=	
		\int_{\pi}f\otimes \left(f\stackrel{r-i-1}{\frown}g\right)\otimes g.
	\end{align}
If both kernels are fully symmetric, this inner product equals
$\tau\left[\left(I_{q+r-2i-2}\left(f\stackrel{i+1}\frown g\right)\right)^2\right]$.
Let $x$ and $s$ each have $r-i-1$ coordinates, let $t$ and $t'$ each have $i+1$ coordinates, and let $y$ have $q-r$ coordinates. Full symmetry gives
\begin{align*}
&\tau\left[\left(I_{q+r-2i-2}\left(f\stackrel{i+1}\frown g\right)\right)^2\right]\\
&=\int f\left(x,y,t\right)g\left(t,s\right)f\left(s,y,t'\right)g\left(t',x\right)\,dx\,dy\,ds\,dt\,dt'\\
&=\int f\left(t,y,x\right)g\left(x,t'\right)f\left(t',y,s\right)g\left(s,t\right)\,dt\,dy\,dt'\,dx\,ds\\
&=\tau\left[\left(I_{q-r+2i+2}\left(f\stackrel{r-i-1}\frown g\right)\right)^2\right].
\end{align*}
Here the coordinate block pairs $\left(x,s\right)$ and $\left(t,t'\right)$ are interchanged, using the invariance of $f$ and $g$ under coordinate permutations. Thus the same inner product can also be written as
	\begin{align}
		&\quad \left\langle\int_{\mathbb{R}_+}\left(g_t^r\stackrel{i}{\frown}\widetilde{f}_t^{k}\right)dt,\int_{\mathbb{R}_+}\left(f_t^q\stackrel{i}{\frown}\widetilde{g}_t^{k-q+r}\right)dt\right\rangle_{L^2\left(\mathbb{R}_+^{q+r-2-2i}\right)}\\
		&=\left\langle f\stackrel{r-i-1}{\frown}g, \left(f\stackrel{r-i-1}{\frown}g\right)^*\right\rangle _{L^2\left(\mathbb{R}_+^{q-r+2i+2}\right)}\notag\\
		&=\tau\left[\left(I_{q-r+2i+2}\left(f\stackrel{r-i-1}{\frown}g\right)\right)^2\right].
	\end{align}
In the second case, combining the kernels gives
	\begin{equation*}
\left[I_{r+k-2i-2}\otimes I_{r-k}\right]\left(f_t^q\stackrel{q-r+i}\frown\widetilde f_t^{q-r+k}+g_t^r\stackrel i\frown\widetilde g_t^k\right).
	\end{equation*} 
Proceeding as above for the two diagonal terms, we obtain a pairing $\pi\in\mathcal P_2\left(q\otimes\left(2i+2\right)\otimes q\right)$ such that
	\begin{equation}
		\left\Vert \int_{\mathbb{R}_+}f_t^q\stackrel{i}{\frown}\widetilde{f}_t^kdt\right\Vert^2_{L^2\left(\mathbb{R}_+^{2q-2i-2}\right)}=\int_{\pi}f\otimes \left(f\stackrel{q-i-1}{\frown}f\right)\otimes f,
	\end{equation}
and a pairing $\pi\in\mathcal P_2\left(r\otimes\left(2i+2\right)\otimes r\right)$ such that
	\begin{equation*}
		\left\Vert \int_{\mathbb{R}_+}g_t^r\stackrel{i}{\frown}\widetilde{g}_t^kdt\right\Vert^2_{L^2\left(\mathbb{R}_+^{2r-2i-2}\right)}=\int_{\pi}g\otimes \left(g\stackrel{r-i-1}{\frown}g\right)\otimes g.
	\end{equation*}
For their inner product, the coordinate expansion is
	
	\begin{align*}
		&\quad \left\langle \int_{\mathbb{R}_+}f_t^q\stackrel{q-r+i}{\frown} \widetilde{f}_t^{q-r+k}dt,\int_{\mathbb{R}_+}g_t^r\stackrel{i}{\frown} \widetilde{g}_t^{k}dt\right\rangle_{L^2\left(\mathbb{R}^{2r-2i-2}_+\right)}\\
		&=\int_{\mathbb R_+^{q+r}} f\left(x_1,\cdots,x_{r-i-1},u_{q-r+i},\cdots,u_{1},t\right)\\
		&\quad \widetilde{f}\left(u_1,\cdots,u_{q-r+i},x_{r-i},\cdots,x_{r+k-2i-2},t,x_{r+k-2i-1},\cdots,x_{2r-2i-2}\right)\\
		&\quad \overline{g\left(x_1,\cdots,x_{r-i-1},u_i',\cdots,u_1',t'\right)}\\
		&\quad \overline{\widetilde{g}\left(u_1',\cdots,u_i',x_{r-i},\cdots,x_{r+k-2i-2},t',x_{r+k-2i-1},\cdots,x_{2r-2i-2}\right)}\\
		&\quad du_1\cdots du_{q-r+i}dt du_1'\cdots du_i'dt'dx_1\cdots dx_{2r-2i-2}\\
		&=\int_{\mathbb R_+^{q+r}} f\left(x_1,\cdots,x_{r-i-1},u_{q-r+i},\cdots,u_{1},t\right)\\
		&\quad \overline{f\left(x_{r+k-2i-2},\cdots,x_{r-i},u_{q-r+i},\cdots,u_1,t,x_{2r-2i-2},\cdots,x_{r+k-2i-1}\right)}\\
		&\quad \overline{g\left(x_1,\cdots,x_{r-i-1},u_i',\cdots,u_1',t'\right)}\\
		&\quad g\left(x_{r+k-2i-2},\cdots,x_{r-i},u_i',\cdots,u_1',t',x_{2r-2i-2},\cdots,x_{r+k-2i-1}\right)\\
		&\quad du_1\cdots du_{q-r+i}dt du_1'\cdots du_i'dt'dx_1\cdots dx_{2r-2i-2}\\
		&=\int_{\mathbb R_+^{q+r}} f\left(x_1,\cdots,x_{r-i-1},u_{q-r+i},\cdots,u_{1},t\right)\\
		&\quad f\left(x_{r+k-2i-1},\cdots,x_{2r-2i-2},t,u_1,\cdots,u_{q-r+i},x_{r-i},\cdots,x_{r+k-2i-2}\right)\\
		&\quad g\left(t',u_1',\cdots,u_i',x_{r-i-1},\cdots,x_1\right)\\
		&\quad g\left(x_{r+k-2i-2},\cdots,x_{r-i},u_i',\cdots,u_1',t',x_{2r-2i-2},\cdots,x_{r+k-2i-1}\right)\\
		&\quad du_1\cdots du_{q-r+i}dt du_1'\cdots du_i'dt'dx_1\cdots dx_{2r-2i-2}.\\			\end{align*}
Integrating with respect to $x_1,\ldots,x_{r-i-1}$ gives a pairing $\pi\in\mathcal P_2\left(q\otimes\left(q-r+2i+2\right)\otimes r\right)$ satisfying
	\begin{align}
		&\quad \left\langle \int_{\mathbb{R}_+}f_t^q\stackrel{q-r+i}{\frown} \widetilde{f}_t^{q-r+k}dt,\int_{\mathbb{R}_+}g_t^r\stackrel{i}{\frown} \widetilde{g}_t^{k}dt\right\rangle_{L^2\left(\mathbb{R}^{2r-2i-2}_+\right)}\notag\\
		&=\int_{\pi} f\otimes \left(g\stackrel{r-i-1}{\frown}f\right)\otimes g.
	\end{align}	
If $f$ and $g$ are fully symmetric, this reduces to
	\begin{align}
		&\quad \left\langle \int_{\mathbb{R}_+}f_t^q\stackrel{q-r+i}{\frown} \widetilde{f}_t^{q-r+k}dt,\int_{\mathbb{R}_+}g_t^r\stackrel{i}{\frown} \widetilde{g}_t^{k}dt\right\rangle_{L^2\left(\mathbb{R}^{2r-2i-2}_+\right)}\notag\\
		&=\left\langle f\stackrel{r-i-1}{\frown} g,f\stackrel{r-i-1}{\frown}g \right\rangle_{L^2\left(\mathbb{R}^{q-r+2i+2}_+\right)}\notag\\
&=\tau\left[I_{q-r+2i+2}\left(f\stackrel{r-i-1}{\frown}g\right)^*I_{q-r+2i+2}\left(f\stackrel{r-i-1}{\frown}g\right)\right].
	\end{align}
Finally, the index pairs $k=q$, $i=q-1$ and $k=r$, $i=r-1$ give the scalar contributions
	\begin{align*}
		\int_{\mathbb{R}_+}\left[I_{q-1+k-1-2i}\otimes I_{q-k}\right]\left(f_t^q\stackrel{i}{\frown}\widetilde{f}_t^k\right)dt&=\left\Vert f\right\Vert^2_{L^2\left(\mathbb{R}_+^q\right)}1_{\mathcal{A}}\otimes 1_{\mathcal{A}}=\tau\left(X^2\right)1_{\mathcal{A}}\otimes 1_{\mathcal{A}},\\
		\int_{\mathbb{R}_+} \left[I_{r-1+k-1-2i}\otimes I_{r-k}\right]\left(g_t^r\stackrel{i}{\frown}\widetilde{g}_t^k\right)dt&=\left\Vert g\right\Vert^2_{L^2\left(\mathbb{R}_+^r\right)}1_{\mathcal{A}}\otimes 1_{\mathcal{A}}=\tau\left(Y^2\right)1_{\mathcal{A}}\otimes 1_{\mathcal{A}}. 
	\end{align*}
	
Therefore, for mirror-symmetric $f$ and $g$, the full expansion is
	\begin{align*}		&\ \ \ \ \ \ \ \ \left\Vert\int_{\mathbb{R}_+}\left(\mbox{id}\otimes\tau\left(\nabla_tF\right)\right)\cdot \left(\nabla_tF\right)^*dt-1_{\mathcal{A}}\otimes 1_{\mathcal{A}}\right\Vert^2_{L^2\left(\mathcal{A},\tau\right)\otimes L^2\left(\mathcal{A},\tau\right)}\\
		&=\sum_{0\le i<k\le q,\ \left(i,k\right)\neq \left(q-1,q\right)}\left\Vert\int_{\mathbb{R}_+}\left(f_t^q\stackrel{i}{\frown}\widetilde{f}_t^k\right)dt\right\Vert_{L^2\left(\mathbb{R}_+^{2q-2i-2}\right)}^2\\
		&\quad +\sum_{k=1}^{q-1}\sum_{i=0}^{\min\left(k-1,r-1\right)}\left\Vert\int_{\mathbb{R}_+}\left(g_t^r\stackrel{i}{\frown}\widetilde{f}_t^k\right)dt\right\Vert_{L^2\left(\mathbb{R}_+^{q+r-2i-2}\right)}^2\\
		&\quad +\sum_{i=0}^{r-1}\left\Vert\int_{\mathbb{R}_+}\left(g_t^r\stackrel{i}{\frown}\widetilde{f}_t^q\right)dt\right\Vert_{L^2\left(\mathbb{R}_+^{q+r-2i-2}\right)}^2\\
		&\quad +\sum_{k=1}^{r-1}\sum_{i=0}^{k-1}\left\Vert\int_{\mathbb{R}^+}\left(f_t^q\stackrel{i}{\frown}\widetilde{g}_t^k\right)dt\right\Vert_{L^2\left(\mathbb{R}_+^{q+r-2i-2}\right)}^2+\sum_{i=0}^{r-1}\left\Vert\int_{\mathbb{R}^+}\left(f_t^q\stackrel{i}{\frown}\widetilde{g}_t^r\right)dt\right\Vert_{L^2\left(\mathbb{R}_+^{q+r-2i-2}\right)}^2\\
		&\quad +\sum_{0\le i<k\le r,\left(i,k\right)\neq \left(r-1,r\right)}\left\Vert\int_{\mathbb{R}_+}\left(g_t^r\stackrel{i}{\frown}\widetilde{g}_t^k\right)dt\right\Vert_{L^2\left(\mathbb{R}_+^{2r-2i-2}\right)}^2\\
		&\quad +2\operatorname{Re}\sum_{k=q-r+1}^{q}\sum_{i=0}^{k-q+r-1}\left\langle\int_{\mathbb{R}_+}\left(g_t^r\stackrel{i}{\frown}\widetilde{f}_t^{k}\right)dt,\int_{\mathbb{R}_+}\left(f_t^q\stackrel{i}{\frown}\widetilde{g}_t^{k-q+r}\right)dt\right\rangle_{L^2\left(\mathbb{R}_+^{q+r-2-2i}\right)}\\
		&\quad +2\operatorname{Re}\sum_{0\le i<k\le r,\left(i,k\right)\neq \left(r-1,r\right)}	\left\langle \int_{\mathbb{R}_+}f_t^q\stackrel{q-r+i}{\frown} \widetilde{f}_t^{q-r+k}dt,\int_{\mathbb{R}_+}g_t^r\stackrel{i}{\frown} \widetilde{g}_t^{k}dt\right\rangle_{L^2\left(\mathbb{R}^{2r-2i-2}_+\right)}\\
	\end{align*}
By the Cauchy--Schwarz inequality and Lemma \ref{estimate of pairing integration}, we have
\begin{align*}
&\sum_{i=0}^{q-2}\sum_{k=i+1}^q\left\Vert\int_0^\infty f_t^q\stackrel i\frown\widetilde f_t^k\,dt\right\Vert_2^2\\
&\le\Vert f\Vert_2^2\sum_{i=0}^{q-2}\left(q-i\right)\left\Vert f\stackrel{q-i-1}\frown f\right\Vert_2
=\Vert f\Vert_2^2\sum_{l=1}^{q-1}\left(l+1\right)\left\Vert f\stackrel l\frown f\right\Vert_2\\
&\le q\Vert f\Vert_2^2\sqrt{q-1}\left(\sum_{l=1}^{q-1}\left\Vert f\stackrel l\frown f\right\Vert_2^2\right)^{\frac{1}{2}}
=q\Vert f\Vert_2^2\sqrt{q-1}\,\kappa_4\left(X\right)^{\frac{1}{2}}.
\end{align*}
For fixed $i$, there are $q-i$ terms, with $k=i+1,\ldots,q$. Only the scalar term $i=q-1,k=q$ is removed; all other terms with $k=q$ are retained.
Similarly,
\begin{align*}
&\sum_{i=0}^{r-2}\sum_{k=i+1}^r\left\Vert\int_0^\infty g_t^r\stackrel i\frown\widetilde g_t^k\,dt\right\Vert_2^2\\
&\le\Vert g\Vert_2^2\sum_{l=1}^{r-1}\left(l+1\right)\left\Vert g\stackrel l\frown g\right\Vert_2
\le r\Vert g\Vert_2^2\sqrt{r-1}\,\kappa_4\left(Y\right)^{\frac{1}{2}}.
\end{align*}
For the mixed integrals $f_t^q\stackrel i\frown\widetilde g_t^k$, the pairing-integral representation obtained above gives
\begin{align*}
&\sum_{i=0}^{r-1}\sum_{k=i+1}^r\left\Vert\int_0^\infty f_t^q\stackrel i\frown\widetilde g_t^k\,dt\right\Vert_2^2\\
&\le\Vert g\Vert_2^2\sum_{i=0}^{r-1}\left(r-i\right)\left\Vert f\stackrel{q-i-1}\frown f\right\Vert_2\\
&\le\Vert g\Vert_2^2\left(\sum_{i=0}^{r-1}\left(r-i\right)^2\right)^{\frac{1}{2}}
\left(\sum_{i=0}^{r-1}\left\Vert f\stackrel{q-i-1}\frown f\right\Vert_2^2\right)^{\frac{1}{2}}
\le r\Vert g\Vert_2^2\sqrt r\,\kappa_4\left(X\right)^{\frac{1}{2}}.
\end{align*}
Here the indices $q-i-1$ are distinct members of $1,\ldots,q-1$.
For the other family of mixed integrals, first consider the terms with $i\le r-2$. The same calculation gives
\begin{align*}
&\sum_{i=0}^{r-2}\sum_{k=i+1}^q\left\Vert\int_0^\infty g_t^r\stackrel i\frown\widetilde f_t^k\,dt\right\Vert_2^2\\
&\le\Vert f\Vert_2^2\sum_{i=0}^{r-2}\left(q-i\right)\left\Vert g\stackrel{r-i-1}\frown g\right\Vert_2\\
&=\Vert f\Vert_2^2\sum_{l=1}^{r-1}\left(q-r+1+l\right)\left\Vert g\stackrel l\frown g\right\Vert_2
\le q\Vert f\Vert_2^2\sqrt{r-1}\,\kappa_4\left(Y\right)^{\frac{1}{2}}.
\end{align*}
When $r=1$, all the sums above whose upper limit is $r-2$ are empty.
It remains to calculate the terms with $i=r-1$. In this case $r-i-1=0$, so the corresponding self-contraction is $g\otimes g$.

Put $k=r+j$, where $0\le j\le q-r$. The coordinate representation of the kernels obtained above gives
\begin{align}
&\left(\int_0^\infty g_t^r\stackrel{r-1}\frown\widetilde f_t^{r+j}\,dt\right)
 \left(y_1,\ldots,y_j,z_1,\ldots,z_{q-r-j}\right)\notag\\
&\qquad=\overline{D_j\left(f,g\right)\left(y_j,\ldots,y_1,z_{q-r-j},\ldots,z_1\right)}.
\label{terminal-internal-identification}
\end{align}
An empty coordinate block is omitted. Reversing the coordinates preserves the $L^2$ norm of the kernel, so the sum of squared norms of these terms is $\sum_{j=0}^{q-r}\Vert D_j\left(f,g\right)\Vert_2^2$.
By the definition of contraction and the mirror symmetry of $g$,
\begin{align*}
&\left(f\stackrel r\frown g\right)\left(x_1,\ldots,x_{q-r}\right)\\
&=\int f\left(x_1,\ldots,x_{q-r},t_r,\ldots,t_1\right)g\left(t_1,\ldots,t_r\right)\,dt_1\cdots dt_r\\
&=\int f\left(x_1,\ldots,x_{q-r},t_1,\ldots,t_r\right)\overline{g\left(t_1,\ldots,t_r\right)}\,dt_1\cdots dt_r\\
&=D_{q-r}\left(f,g\right)\left(x_1,\ldots,x_{q-r}\right).
\end{align*}
On the other hand, using the mirror symmetry of $f$ and reversing $t_1,\ldots,t_r$, we obtain
\begin{align*}
&\left(f\stackrel r\frown g\right)^*\left(x_1,\ldots,x_{q-r}\right)\\
&=\int\overline{f\left(x_{q-r},\ldots,x_1,t_1,\ldots,t_r\right)}g\left(t_1,\ldots,t_r\right)\,dt_1\cdots dt_r\\
&=\int f\left(t_r,\ldots,t_1,x_1,\ldots,x_{q-r}\right)g\left(t_1,\ldots,t_r\right)\,dt_1\cdots dt_r\\
&=\int f\left(t_1,\ldots,t_r,x_1,\ldots,x_{q-r}\right)\overline{g\left(t_1,\ldots,t_r\right)}\,dt_1\cdots dt_r\\
&=D_0\left(f,g\right)\left(x_1,\ldots,x_{q-r}\right).
\end{align*}
Thus the two endpoint terms have equal squared norms, and
\begin{align*}
&\sum_{i=0}^{r-1}\sum_{k=i+1}^q\left\Vert\int_0^\infty g_t^r\stackrel i\frown\widetilde f_t^k\,dt\right\Vert_2^2\\
&=\sum_{i=0}^{r-2}\sum_{k=i+1}^q\left\Vert\int_0^\infty g_t^r\stackrel i\frown\widetilde f_t^k\,dt\right\Vert_2^2
 +2\left\Vert f\stackrel r\frown g\right\Vert_2^2+\mathcal R_{q,r}\left(f,g\right).
\end{align*}
To calculate the endpoint norm, let $x$ be a block of $q-r$ integration variables, and let $t,s$ each be a block of $r$ integration variables. Then
\begin{align}
\left\Vert f\stackrel r\frown g\right\Vert_2^2
&=\int\left|\int f\left(x,t\right)\overline{g\left(t\right)}\,dt\right|^2\,dx\notag\\
&=\int f\left(x,t\right)\overline{f\left(x,s\right)}\overline{g\left(t\right)}g\left(s\right)\,dx\,dt\,ds\notag\\
&=\int\left(f\stackrel{q-r}\frown f\right)\left(s_r,\ldots,s_1,t_1,\ldots,t_r\right)
 \overline{g\left(t_1,\ldots,t_r\right)}g\left(s_1,\ldots,s_r\right)\,dt\,ds\notag\\
&=\left\langle f\stackrel{q-r}\frown f,g\otimes g\right\rangle_2\notag\\
&\le\Vert g\Vert_2^2\left\Vert f\stackrel{q-r}\frown f\right\Vert_2
\le\Vert g\Vert_2^2\,\kappa_4\left(X\right)^{\frac{1}{2}}.
\label{endpoint-contraction-bound}
\end{align}
The third equality uses $\overline{f\left(x,s\right)}=f\left(s_r,\ldots,s_1,x_{q-r},\ldots,x_1\right)$.
For the fourth equality, we reverse $s_1,\ldots,s_r$ and use the mirror symmetry of $g$.
The last inequality follows from $1\le q-r\le q-1$.

For the two families of cross inner products in the expansion above, each ordered pair of chaos degrees occurs at most twice, either among the two diagonal expansions or among the two mixed expansions. Thus each term is combined with at most one term from the other expansion.
Observe that the tensor-valued Wigner integrals corresponding to $1\le j\le q-r-1$ are
\[
\left[I_j^S\otimes I_{q-r-j}^S\right]
\left(\int_0^\infty g_t^r\stackrel{r-1}\frown\widetilde f_t^{r+j}\,dt\right).
\]
These terms have distinct ordered pairs of chaos degrees, and the sum of their squared norms is exactly $\mathcal R_{q,r}\left(f,g\right)$.
Their total chaos degree is the odd integer $q-r$, and both components have positive degree. They therefore cannot coincide with diagonal terms.
The mixed terms with $i\le r-2$ have total degree at least $q-r+2$. In the other mixed expansion, the only term of total degree $q-r$ has ordered degrees $\left(q-r,0\right)$.
The two endpoint terms have ordered degrees $\left(0,q-r\right)$ and $\left(q-r,0\right)$, respectively.
Hence these strictly internal terms are orthogonal to every remaining term.
Combining all the estimates above, we obtain
\begin{align}
\Sigma^*\left(F\mid S\right)^2
&\le\left\Vert\int_0^\infty
 \left[\left(\operatorname{id}\otimes\tau\right)\left(\nabla_tF\right)\otimes1\right]\left(\nabla_tF\right)^*\,dt-1\otimes1\right\Vert_2^2\notag\\
&\le\mathcal R_{q,r}\left(f,g\right)
 +2\left[q\Vert f\Vert_2^2\sqrt{q-1}+\Vert g\Vert_2^2\left(r\sqrt r+2\right)\right]\kappa_4\left(X\right)^{\frac{1}{2}}\notag\\
&\quad+2\left(r\Vert g\Vert_2^2+q\Vert f\Vert_2^2\right)\sqrt{r-1}\,\kappa_4\left(Y\right)^{\frac{1}{2}}\notag\\
&\le\mathcal R_{q,r}\left(f,g\right)+2A_{q,r}\kappa_4\left(X\right)^{\frac{1}{2}}
 +2q\sqrt{r-1}\,\kappa_4\left(Y\right)^{\frac{1}{2}}\notag\\
&\le\mathcal R_{q,r}\left(f,g\right)
 +2\sqrt{A_{q,r}^2+q^2\left(r-1\right)}\left(\kappa_4\left(X\right)+\kappa_4\left(Y\right)\right)^{\frac{1}{2}}\notag\\
&\le\mathcal R_{q,r}\left(f,g\right)+C_{q,r}\kappa_4\left(F\right)^{\frac{1}{2}}.
\label{proof-remainder-discrepancy-bound}
\end{align}
Here we use $\Vert f\Vert_2^2+\Vert g\Vert_2^2=1$, $r\le q$, and the inequality $\kappa_4\left(X\right)+\kappa_4\left(Y\right)\le\kappa_4\left(F\right)$ proved above.
Together with $W_2\left(F,S\right)\le\Sigma^*\left(F\mid S\right)$, this proves the bound with a remainder under mirror symmetry.

If the higher-order kernel $f$ is also cyclically symmetric, moving its first $j+r$ coordinates cyclically to the end gives
\begin{align}
&D_j\left(f,g\right)\left(x_1,\ldots,x_{q-r}\right)\notag\\
&=\int f\left(x_{j+1},\ldots,x_{q-r},x_1,\ldots,x_j,t_1,\ldots,t_r\right)
 \overline{g\left(t_1,\ldots,t_r\right)}\,dt_1\cdots dt_r\notag\\
&=\left(f\stackrel r\frown g\right)\left(x_{j+1},\ldots,x_{q-r},x_1,\ldots,x_j\right).
\label{cyclic-internal-permutation}
\end{align}
This change merely permutes the free coordinates and preserves the squared norm. Thus all $q-r+1$ terms with $i=r-1$ are controlled by the endpoint norm, and
\begin{align}
&\sum_{k=r}^q\left\Vert\int_0^\infty g_t^r\stackrel{r-1}\frown\widetilde f_t^k\,dt\right\Vert_2^2\notag\\
&=\left(q-r+1\right)\left\Vert f\stackrel r\frown g\right\Vert_2^2
\le\left(q-r+1\right)\Vert g\Vert_2^2\,\kappa_4\left(X\right)^{\frac{1}{2}}.
\label{cyclic-terminal-sum}
\end{align}
Bounding the two families of cross inner products together with the four sums of squared norms now yields
\begin{align*}
\Sigma^*\left(F\mid S\right)^2
&\le2\left[q\Vert f\Vert_2^2\sqrt{q-1}
 +\Vert g\Vert_2^2\left(r\sqrt r+q-r+1\right)\right]\kappa_4\left(X\right)^{\frac{1}{2}}\\
&\quad+2\left(r\Vert g\Vert_2^2+q\Vert f\Vert_2^2\right)\sqrt{r-1}\,\kappa_4\left(Y\right)^{\frac{1}{2}}\\
&\le2A_{q,r}^{\mathrm{cyc}}\kappa_4\left(X\right)^{\frac{1}{2}}
 +2q\sqrt{r-1}\,\kappa_4\left(Y\right)^{\frac{1}{2}}\\
&\le2\sqrt{\left(A_{q,r}^{\mathrm{cyc}}\right)^2+q^2\left(r-1\right)}
 \left(\kappa_4\left(X\right)+\kappa_4\left(Y\right)\right)^{\frac{1}{2}}\\
&\le C_{q,r}^{\mathrm{cyc}}\kappa_4\left(F\right)^{\frac{1}{2}}.
\end{align*}
Applying the free WS inequality again and taking square roots proves the quantitative bound without a remainder under cyclic symmetry.

Now suppose that both $f$ and $g$ are fully symmetric. By the integral calculations above for fully symmetric kernels, for every admissible $i,k$,
\begin{align*}
\left\Vert\int_0^\infty f_t^q\stackrel i\frown\widetilde g_t^k\,dt\right\Vert_2^2
&=\left\Vert f\stackrel{i+1}\frown g\right\Vert_2^2,\\
\left\Vert\int_0^\infty g_t^r\stackrel i\frown\widetilde f_t^k\,dt\right\Vert_2^2
&=\left\Vert g\stackrel{i+1}\frown f\right\Vert_2^2
=\left\Vert f\stackrel{i+1}\frown g\right\Vert_2^2.
\end{align*}
The analogous identities hold for the diagonal terms. Setting $l=i+1$ and counting terms over the same summation ranges as above, we obtain
\begin{align*}
\sum_{i=0}^{q-2}\sum_{k=i+1}^q\left\Vert\int_0^\infty f_t^q\stackrel i\frown\widetilde f_t^k\,dt\right\Vert_2^2&=\sum_{l=1}^{q-1}\left(q-l+1\right)\left\Vert f\stackrel l\frown f\right\Vert_2^2\le q\kappa_4\left(X\right),\\
\sum_{i=0}^{r-2}\sum_{k=i+1}^r\left\Vert\int_0^\infty g_t^r\stackrel i\frown\widetilde g_t^k\,dt\right\Vert_2^2
&=\sum_{l=1}^{r-1}\left(r-l+1\right)\left\Vert g\stackrel l\frown g\right\Vert_2^2\le r\kappa_4\left(Y\right).
\end{align*}
The two mixed families contain, respectively, $r-l+1$ and $q-l+1$ terms with the same norm. Hence
\begin{align*}
\sum_{i=0}^{r-1}\sum_{k=i+1}^r\left\Vert\int_0^\infty f_t^q\stackrel i\frown\widetilde g_t^k\,dt\right\Vert_2^2+\sum_{i=0}^{r-1}\sum_{k=i+1}^q\left\Vert\int_0^\infty g_t^r\stackrel i\frown\widetilde f_t^k\,dt\right\Vert_2^2&=\sum_{l=1}^r\left(q+r-2l+2\right)\left\Vert f\stackrel l\frown g\right\Vert_2^2\\
&\le\left(q+r\right)\sum_{l=1}^r\left\Vert f\stackrel l\frown g\right\Vert_2^2.
\end{align*}
Since all the free coordinates of $f\stackrel r\frown g$ belong to $f$, this kernel is real-valued and fully symmetric. Consequently,
$\left(f\stackrel r\frown g\right)^*=f\stackrel r\frown g$.
Substitution into the positive decomposition of the fourth cumulant established above gives
\begin{align*}
\kappa_4\left(F\right)-\kappa_4\left(X\right)-\kappa_4\left(Y\right)
&\ge2\sum_{l=1}^{r-1}\left\Vert f\stackrel l\frown g\right\Vert_2^2
 +8\left\Vert f\stackrel r\frown g\right\Vert_2^2\\
&\ge2\sum_{l=1}^r\left\Vert f\stackrel l\frown g\right\Vert_2^2.
\end{align*}
Combining the estimates for the four sums of squared norms and the cross inner products, we conclude that
\begin{align*}
W_2\left(F,S\right)^2&\le\Sigma^*\left(F\mid S\right)^2\\
&\le2\left[q\kappa_4\left(X\right)+r\kappa_4\left(Y\right)
 +\left(q+r\right)\sum_{l=1}^r\left\Vert f\stackrel l\frown g\right\Vert_2^2\right]\\
&\le2q\left[\kappa_4\left(X\right)+\kappa_4\left(Y\right)
 +2\sum_{l=1}^r\left\Vert f\stackrel l\frown g\right\Vert_2^2\right]
\le2q\kappa_4\left(F\right).
\end{align*}
Taking square roots gives $W_2\left(F,S\right)\le\sqrt{2q}\,\kappa_4\left(F\right)^{\frac{1}{2}}$.
The asserted $W_2$ convergence follows directly from these estimates. Convergence under every fixed $W_p$ follows from the uniform support bound and the independent interpolation argument in Proposition \ref{support-Wp-interpolation} below.
\end{proof}

\begin{remark}
The remainder in \eqref{remainder-definition} contains the strictly internal contractions. The endpoint contractions are controlled by \eqref{endpoint-contraction-bound}. In particular, $q=r+1$ implies $\mathcal R_{q,r}\left(f,g\right)=0$. 
\end{remark}
\begin{remark}
Under cyclic symmetry, \eqref{cyclic-internal-permutation} gives
$D_0\left(f,g\right)=D_{q-r}\left(f,g\right)=f\stackrel r\frown g$.
Together with mirror symmetry, this implies
$\left(f\stackrel r\frown g\right)^*=f\stackrel r\frown g$.
The positive decomposition \eqref{positive-mixed-decomposition} therefore yields
\[
\mathcal R_{q,r}\left(f,g\right)
=\left(q-r-1\right)\left\Vert f\stackrel r\frown g\right\Vert_2^2
\le\frac{q-r-1}{8}\kappa_4\left(F\right).
\]
Thus cyclic symmetry also gives a direct quantitative control of the remainder.
\end{remark}

\section{From \texorpdfstring{$W_2$}{W2} to \texorpdfstring{$W_p$}{Wp} under a support bound}
	For probability measures $\mu$ and $\nu$ on $\mathbb{R}$, let $W_p(\mu,\nu)$ denote the usual $p$-Wasserstein distance. In one variable, the free $W_2$ distance used above agrees with the classical $W_2$ distance between the spectral distributions \cite{BV01}. The following elementary interpolation estimate is the step that permits the passage from $W_2$ to $W_p$.

	\begin{proposition}\label{support-Wp-interpolation}
		Let $\mu$ and $\nu$ be probability measures supported, respectively, in $[-R_\mu,R_\mu]$ and $[-R_\nu,R_\nu]$, and set $D=R_\mu+R_\nu$. Then
		\begin{equation*}
			W_p(\mu,\nu)\le W_2(\mu,\nu),\qquad 1\le p\le2,
		\end{equation*}
		whereas
		\begin{equation*}
			W_p(\mu,\nu)\le D^{1-\frac{2}{p}}W_2(\mu,\nu)^{\frac{2}{p}},\qquad 2<p<\infty.
		\end{equation*}
	\end{proposition}
	\begin{proof}
		Let $\pi$ be any coupling of $\mu$ and $\nu$. If $1\le p\le2$, the monotonicity of $L^p(\pi)$-norms gives
		\begin{equation*}
			\left(\int_{\mathbb{R}^2}|x-y|^p\,d\pi(x,y)\right)^{\frac{1}{p}}
			\le\left(\int_{\mathbb{R}^2}|x-y|^2\,d\pi(x,y)\right)^{\frac{1}{2}}.
		\end{equation*}
		Taking the infimum over $\pi$ proves the first assertion. If $p>2$, then $|x-y|\le D$ on the support of every coupling, and hence
		\begin{equation*}
			\int_{\mathbb{R}^2}|x-y|^p\,d\pi(x,y)
			\le D^{p-2}\int_{\mathbb{R}^2}|x-y|^2\,d\pi(x,y).
		\end{equation*}
		Applying this inequality to an optimal $W_2$-coupling and taking the $p$-th root proves the second assertion.
	\end{proof}

	We now apply Proposition \ref{support-Wp-interpolation} with the same chaos orders $q>r$ as in Theorem \ref{Maintheorem2}. The operator-norm estimate for Wigner integrals \cite{BIA98} gives
\[
\left\Vert I_k^S\left(h\right)\right\Vert_\infty\le\left(k+1\right)\Vert h\Vert_2.
\]
Consequently,
\begin{align*}
\Vert F\Vert_\infty
&\le\left\Vert I_q^S\left(f\right)\right\Vert_\infty
 +\left\Vert I_r^S\left(g\right)\right\Vert_\infty\\
&\le\left(q+1\right)\Vert f\Vert_2+\left(r+1\right)\Vert g\Vert_2.
\end{align*}
The standard semicircular law is supported in $\left[-2,2\right]$. Thus the two support radii have sum at most
\[
D_{q,r}\left(f,g\right)
=\left(q+1\right)\Vert f\Vert_2+\left(r+1\right)\Vert g\Vert_2+2
\le\sqrt{\left(q+1\right)^2+\left(r+1\right)^2}+2,
\]
where the last inequality follows from $\Vert f\Vert_2^2+\Vert g\Vert_2^2=1$.

\begin{corollary}\label{Wp-fourth-moment-corollary}
Under the hypotheses of Theorem \ref{Maintheorem2}, for every $1\le p<\infty$,
\[
W_p\left(F,S\right)\le
\begin{cases}
\left[C_{q,r}\kappa_4\left(F\right)^{\frac{1}{2}}+\mathcal R_{q,r}\left(f,g\right)\right]^{\frac{1}{2}},&1\le p\le2,\\[4pt]
D_{q,r}\left(f,g\right)^{1-\frac{2}{p}}
\left[C_{q,r}\kappa_4\left(F\right)^{\frac{1}{2}}+\mathcal R_{q,r}\left(f,g\right)\right]^{\frac{1}{p}},&2<p<\infty.
\end{cases}
\]
If the higher-order kernel $f$ is also cyclically symmetric, then
\[
W_p\left(F,S\right)\le
\begin{cases}
\left(C_{q,r}^{\mathrm{cyc}}\right)^{\frac{1}{2}}\kappa_4\left(F\right)^{\frac{1}{4}},&1\le p\le2,\\[4pt]
D_{q,r}\left(f,g\right)^{1-\frac{2}{p}}
\left(C_{q,r}^{\mathrm{cyc}}\right)^{\frac{1}{p}}\kappa_4\left(F\right)^{\frac{1}{2p}},&2<p<\infty.
\end{cases}
\]
If $f,g$ are fully symmetric, then
\[
W_p\left(F,S\right)\le
\begin{cases}
\sqrt{2q}\,\kappa_4\left(F\right)^{\frac{1}{2}},&1\le p\le2,\\[4pt]
D_{q,r}\left(f,g\right)^{1-\frac{2}{p}}\left(2q\right)^{\frac{1}{p}}\kappa_4\left(F\right)^{\frac{1}{p}},&2<p<\infty.
\end{cases}
\]
\end{corollary}
\begin{proof}
For $1\le p\le2$, Proposition \ref{support-Wp-interpolation} gives
$W_p\left(F,S\right)\le W_2\left(F,S\right)$, so the corresponding bounds follow from Theorem \ref{Maintheorem2} and Corollary \ref{cyclic-quantitative-corollary}.
For $2<p<\infty$, the same proposition gives
\[
W_p\left(F,S\right)^p
\le D_{q,r}\left(f,g\right)^{p-2}W_2\left(F,S\right)^2.
\]
Substituting the three $W_2$ estimates in turn, we obtain
\begin{align*}
W_p\left(F,S\right)^p
&\le D_{q,r}\left(f,g\right)^{p-2}
 \left[C_{q,r}\kappa_4\left(F\right)^{\frac{1}{2}}+\mathcal R_{q,r}\left(f,g\right)\right]
\end{align*}
in the mirror symmetric case,
\[
W_p\left(F,S\right)^p
\le D_{q,r}\left(f,g\right)^{p-2}C_{q,r}^{\mathrm{cyc}}\kappa_4\left(F\right)^{\frac{1}{2}}
\]
in the cyclic case, and
\[
W_p\left(F,S\right)^p
\le2qD_{q,r}\left(f,g\right)^{p-2}\kappa_4\left(F\right)
\]
in the fully symmetric case. Taking the $p$-th root proves the assertions.
\end{proof}

\begin{corollary}\label{vanishing-remainder-corollary}
Fix positive integers $q>r$ of opposite parity. Let $f_N,g_N$ be mirror symmetric with
$\Vert f_N\Vert_2^2+\Vert g_N\Vert_2^2=1$, and put
$F_N=I_q^S\left(f_N\right)+I_r^S\left(g_N\right)$. Assume
\[
\mathcal R_{q,r}\left(f_N,g_N\right)\longrightarrow0.
\]
Then the following statements are equivalent:
\begin{enumerate}
\item $\kappa_4\left(F_N\right)\to0$, equivalently $\tau\left(F_N^4\right)\to2$;
\item $F_N$ converges in distribution to a standard semicircular variable $S$;
\item $W_p\left(F_N,S\right)\to0$ for every fixed $1\le p<\infty$.
\end{enumerate}
\end{corollary}
\begin{proof}
The operator-norm estimate above gives a common compact interval containing the spectral supports of all $F_N$ and $S$. If the first assertion holds, the quantities
$C_{q,r}\kappa_4\left(F_N\right)^{\frac{1}{2}}+\mathcal R_{q,r}\left(f_N,g_N\right)$ tend to zero, whereas
$D_{q,r}\left(f_N,g_N\right)$ is uniformly bounded. Corollary \ref{Wp-fourth-moment-corollary} therefore implies the third assertion. The third assertion implies convergence in distribution. Finally, weak convergence on a common compact interval implies convergence of the fourth moments, which gives the first assertion.
\end{proof}

	At the end of this section, we provide an example to illustrate that the two chaos orders $q,r$ must have different parity.
	\begin{example}
		Consider
		$f=1_{[0,1]^3}(x,y,z),g=1_{[0,1]}(x)$. Then by Theorem 5.3.4 of \cite{BIA98}, we know that
		\begin{equation*}
			I_3^S(f)=S_1^3-2S_1,\ I_1^S(g)=S_1.
		\end{equation*}
		For real coefficients $\alpha,\beta$, put $Z=\alpha I_3^S(f)+\beta I_1^S(g)$. On the one hand, let
		\begin{equation*}
			\tau(Z^2)=\alpha^2+\beta^2=1.
		\end{equation*}
		On the other hand, by calculation we obtain
		\begin{align*}
			\tau(Z^4)&=\tau\left[\left(\alpha S_1^3-(2\alpha-\beta)S_1\right)^4\right]=4\alpha^4+8\alpha^3\beta+12\alpha^2\beta^2+4\alpha\beta^3+2\beta^4.
		\end{align*}
		Let
		\begin{equation*}
			4\alpha^4+8\alpha^3\beta+12\alpha^2\beta^2+4\alpha\beta^3+2\beta^4=2.
		\end{equation*}
		To obtain a nontrivial example, we seek coefficients with $\alpha\beta\ne0$. Put $t=\dfrac{\alpha}{\beta}$. Using $\alpha^2+\beta^2=1$, the fourth-moment condition is equivalent to
		\begin{equation*}
			2t\left(t^3+4t^2+4t+2\right)=0.
		\end{equation*}
		Since $t\ne0$, this reduces to
		\begin{equation*}
			t^3+4t^2+4t+2=0.
		\end{equation*}
		Choose a root $t\in(-3,-2)$, which exists by the intermediate value theorem, and put $\beta=(1+t^2)^{-\frac{1}{2}}$ and $\alpha=t\beta$. Then
		\begin{equation*}
			\tau(Z^2)=1,\ \tau(Z^4)=2.
		\end{equation*}
		but $Z$ does not have the standard semicircular distribution.
	\end{example}

\end{document}